\documentclass[12pt]{amsart}
\usepackage{amsfonts,amsthm,amsmath,amssymb,multicol,pstcol,pst-plot,
pst-grad,pst-node,graphicx}
\newtheorem{thm}{Theorem}[section]

\newtheorem{ex}[thm]{Example}
\newtheorem{cor}[thm]{Corollary}
\newtheorem{lem}[thm]{Lemma}
\newtheorem{prop}[thm]{Proposition}
\theoremstyle{definition}
\newtheorem{defn}[thm]{Definition}
\theoremstyle{remark}
\newtheorem{rem}[thm]{Remark}
\newtheorem{notat}[thm]{Notation}
\numberwithin{equation}{section}

\newcommand{\A}{\mathcal{A}}
\newcommand{\comment}[1]{}

\def\PNC{\mathcal{NCF}}
\def\T{\textrm{T-}}
\def\Tr{\mbox{Tr}\;}

\def\l{\lambda}
\def\P{\mathbf{P}}

\def\tr{\mbox{Tr} \;}

\def\Frtn1{ {{\frac{1}{2^{n+1}}}} }
\newfont{\scal}{cmsy10 scaled600}
\newfont{\pic}{cmb10 scaled950}

\newfont{\CAL}{cmsy10 scaled1700}
\newfont{\Cal}{cmsy10 scaled850}
\newcommand{\iso}{\simeq}
\newcommand{\isosup}{\supset_{\simeq}}

\newcommand{\isosub}{\subset_{\simeq}}

\comment{
\definecolor{Mygray}{gray}{.92}
\definecolor{Mygray1}{gray}{.88}
\definecolor{Mygray2}{gray}{.85}
\definecolor{Mygray3}{gray}{.96}
\newfont{\scal}{cmsy10 scaled600}
\newfont{\pic}{cmb10 scaled950}
\newfont{\CAL}{cmsy10 scaled1700}
\newfont{\Cal}{cmsy10 scaled850}
}

\newcounter{alphnum}

\newcounter{Romnum}

\begin{document}

\date{December 12, 2006 }
\title[]{The Structure of Two-parabolic Space:\\ Parabolic Dust and
Iteration.}

\author{Jane Gilman}

\thanks{Research supported in part by grants from the NSA and the Rutgers
Research Council and by Yale University.}



%


\begin{abstract}
A non-elementary M\"obius group generated by two-parabolics is
determined up to conjugation by one complex parameter and the
parameter space has been extensively studied. In this paper, we use
the results of \cite{GW} to obtain an additional structure for the
parameter space, which we term the {\sl two-parabolic space}. This
structure allows us to identify groups that contain additional
conjugacy classes of primitive parabolics, which following
\cite{Indra} we call {\sl parabolic dust groups}, non-free groups
off the real axis, and groups that are both parabolic dust and
non-free; some of these contain $\mathbb{Z} \times \mathbb{Z}$
subgroups. The structure theorem also attaches additional geometric
structure to discrete and non-discrete groups lying in given regions
of the parameter space including a new explicit construction of some
non-classical $\T$Schottky groups.
\end{abstract}

\maketitle

\section{Introduction} Non-elementary M\"obius groups generated by two parabolic
transformations are parameterized by a non-zero complex number,
$\lambda$. This parameter space has been studied extensively. The
study begins with the work of Lyndon and Ullman \cite{LU} who
obtained results about the {\sl free part}, those values of $\lambda
\in \mathbb{C}$ corresponding to free groups and David Wright who,
in the 1970's,  obtained computer pictures of the portion of the
discrete free part that has come to be known as the Riley slice,
$\mathcal{R}$ \cite{Indra, Wr}  (see also \cite{Chang, Ree}).
Subsequently Keen and Series studied the Riley slice and its
boundary \cite{KS}.  Beardon has extensive results on free and
non-free points \cite{BPell, Beard}. Bamberg \cite{BAM} also studied
the non-free points. Gehring, Machlachlan, and Martin \cite{GMM}
have studied points for which $\lambda$ gives an arithmetic group.
Recently Agol \cite{Ag} has
classified the non-free discrete groups via the topology of a
corresponding knot. The discrete two parabolic groups are, of
course, Kleinian groups.

 In \cite{GW}, Gilman and Waterman found
explicit equations for the boundary for the subspace consisting of
classical T-Schottky groups. The boundary consisted of the portions
of two parabolas between their intersection. \comment{A similar
boundary was found for the space of non-separating disjoint circle
groups, NSDC groups for short.} \comment{This boundary consisted of
portions of two parabolas. and a similar boundary for the subspace
of the classical T-Schottky groups.} In addition they found a one
complex parameter family of non-classical T-Schottky groups. NSDC
groups, classical T-Schottky groups and non-classical T-Schottky
groups are all either inside the Riley slice or on its boundary.
Groups interior to the Riley slice correspond to discrete groups
where there are two, but not three conjugacy classes of maximal
parabolic subgroups, that is, quotient surfaces or handlebodies
where two, but not three curves have been {\sl pinched}. In
\cite{GW} it was shown that only four values of $\l$ on the boundary
of the classical $\T$Schottky space correspond to additional
pinching. These are $\l = \pm i, \pm 2$.

 Here we are interested in all groups whether or not they are
discrete and we want to identify classes of non-discrete groups by
the algebraic properties of $\l$. These algebraic properties are
strongly related to geometric properties of the action of the group
on $\hat{\mathbb{C}}$.

We use the results of \cite{GW} to obtain an additional structure
theorem for the parameter space, which we term the {\sl
two-parabolic space}. This structure allows us to identify groups
that contain additional  conjugacy classes of primitive parabolics,
which following \cite{Indra} we call parabolic dust groups, non-free
groups off the real axis,  and groups that are simultaneously
parabolic dust and non-free.


\section{Organization}
The organization of this paper is as follows. Section
\ref{section:prelimnot} contains some background, notation and
terminology. Some readers may skip or skim this section. The body of
the paper begins with section \ref{section:par} where the relation
between $\l$ and the matrices of the generators is defined and
iteration of $\l$ is introduced. Section \ref{section:schottprior}
reviews definitions for Schottky groups and prior results about the
classical $\T$Schottky boundary. The results of these two sections
are combined to obtain the idea of $n$th classical $\T$Schottky
groups and the main results about such groups (Theorems
\ref{theorem:tncircle},\ref{theorem:iterates},
\ref{theorem:whitehead} and \ref{theorem:nonfree} and corollaries
\ref{cor:or} and  \ref{cor:generalwhite}). Section
\ref{section:tess} describes the tessellation of the complex plane
by the classical $\T$Schottky boundary and its pre-images and gives
some example. Section \ref{section:lateral} introduces a new concept
of iteration, lateral iteration and discusses moving around the
parameter space using combined lateral and vertical iteration
(proposition \ref{prop:lv}). The regions for which we find equations
and boundaries correspond to types of configuration of Schottky
curves on $\hat{\mathbb{C}}$, that is to the geometry of the action
of the generators of $G_\l$ (section
\ref{section:configs}). 
We close with some ideas about generalizing the lateral and vertical
moves to arbitrary two-generator groups (section
\ref{section:conclude}).

\tableofcontents
\section{Preliminaries: Background, terminology and notation}
\label{section:prelimnot}

We let $\mathbb{M}$ be the M{\"o}bius group, the group of two-by-two
nonsingular complex matrices acting as linear transformations on the
complex sphere $\hat{\mathbb{C}}$ and by extension on hyperbolic
three space. \comment{ $\mathbb{H}^3$. We identify $\mathbb{H}^3$
with upper-half-three space $\{(x,y,t) \in {\mathbb{R}}^3 | t>0\}$.
The hyperbolic metric on $\mathbb{H}^3$ can be defined as the unique
metric that makes the M{\"o}bius group the full group of orientation
preserving isometries \cite{Wilk}. }
We identify $\mathbb{M}$ with
$PSL(2, \mathbb{C})$.

We are interested in two generator subgroups,  $G_\l$,  where $\l$
is a non-zero complex number, {\sl  the parameter}, and $G_\l$ is
generated by two
 parabolic transformations.  \comment{ and often work with the preimage of
a M\"obius transformation in $SL(2,\mathbb{C})$.}  Although the
trace of an element of $PSL(2,\mathbb{C})$
 is only determined up to sign, if we require that the pull backs to $SL(2,\mathbb{C})$
 of the two generators of $G_\l$ have positive trace, the trace of every other
 element of $G_\l$ is well-defined. We use $\tr$ to denote the trace
 of an element of $G_\l$ or its appropriate pull-back to an element
 of $SL(2,\mathbb{C})$. A transformation is parabolic, of course,
 if its trace is two.
A group is {\sl non-elementary} if it contains no abelian subgroups
of finite index. The precise matrices that generate $G_\l$ will be
given below (section \ref{section:par}).

\subsection{The sometimes topological picture} We let $\Omega(G_\l)$ denote set
of discontinuity of $G_\l$. The group $G_\l$ may or may not be
discrete, but if it is discrete it can have at most three conjugacy
classes of maximal parabolic
subgroups.  
The {\sl Riley slice of Schottky space} is denoted by $\mathcal{R}$
and is defined to be the set of $\l$ where $ G_\l$ is discrete and
free and contains exactly two conjugacy classes of primitive
parabolics  so that $\Omega(G_\l)/G_\l$ is a four punctured sphere.
We note that when
such a $G_\l$ is discrete,  
the corresponding hyperbolic three manifold will be doubly cusped at
the images of the parabolic fixed points of the group and the
corresponding boundary of the handlebody will be a four punctured
sphere where the punctures are identified in pairs, that is,  a
surface of genus two where two curves have been {\sl pinched} to a
point.

When the group $G_\l$ is not discrete, there is no topological
picture to think of.

\subsection{Properties} We will be interested in describing sets of $\l$
where $G_\l$ has a number of different algebraic and/or geometric
properties. By abuse of language, we say that $\l$ is free meaning
that $G_\l$ is a free group. We say that $\l$ is Schottky,
T-Schottky, classical T-Schottky, non-classical, non-free,  NSDC,
etc.  if $G_\l$ is.

\subsection{Terminology}

If $G$ is any subgroup of $\mathbb{M}$ we call  a parabolic element
$h \in G$  a {\sl primitive parabolic}  if $h = g^m$ for some $g \in
G$ and some integer $m$ implies $m = \pm 1$. If $G$ is a Kleinian
group acting on $\hat{\mathbb{C}}$ and $p$ is any point fixed by an
element of $G$, then $G_p$, the stabilizer of $p$ is either a rank 1
or rank 2 parabolic subgroup. A Kleinian group is {\sl maximally
parabolic} if it allows no deformations with more conjugacy classes
of parabolic elements.

The group $G_\l$ might not be discrete and have three or more
conjugacy classes of primitive parabolics. Following \cite{Indra}
(page 258) and Wright's 2004 post-card caption, we term a two
parabolic generator group, whether it is discrete or not, parabolic
{\sl dust}  if it has a third primitive conjugacy class of
parabolics. In that case both the group and its parameter $\l$ are
referred to as {\sl dust}.

\subsection{Parabolics and pinching}
If $G$ is a discrete group so that there exists a quotient surface,
$S= \Omega(G)/G$  and a quotient manifold,  words in $G$ correspond
curves on $S$. If $G$ has a deformation sending some non-parabolic
word  to a parabolic, we say that the corresponding curve has been
{\sl pinched}.

 The
 maximal number of elements that can be pinched to parabolics is
 precisely the number of rank $1$ parabolic subgroups that any
 group isomorphic to $G$ can contain. A group with this largest
 number of rank $1$ parabolic subgroups is maximally
 parabolic. If $G$ is a Schottky group of rank two, then
it can contain at most three rank $1$ parabolic subgroups.

\section{The matrices and iterating $\l$} \label{section:par}

If $G$ is a marked group with two  parabolic generators, $S$ and
$T$, then up to conjugation $G =G_\lambda = \langle S,T \rangle $
where $S=\left( \begin{array}{cc}
1 & 0\\
1 & 1
\end{array}
\right)$
 and
$T=T_\lambda = \left( \begin{array}{cc}
1 & 2\lambda\\
0 & 1
\end{array}
\right)$ for some nonzero complex number.
We assume that $\lambda \ne 0$ so that $G$ is non-elementary.

We note for future reference that since $\tr \; S =2, \tr \; T =2$
and $\lambda \ne 0$, letting $[S,T]$ denote the multiplicative
commutator, we have
\begin{equation} \label{equation:Glam}
\tr \; [S,T] -2= \tr \;  STS^{-1}T^{-1} -2 = 4 \lambda^2;
\end{equation}
\begin{equation}
\tr \;  ST^{-1} = \tr \; TS^{-1} = 2 - 2\lambda \end{equation}
\begin{equation} \mbox{ and }\\ \;\;\tr \; ST = \tr \; (TS)^{-1} = 2 +
2\lambda .
\end{equation}
\begin{equation}
 \tr\; ST =\pm 2 \Leftrightarrow
\lambda = -2, \;  \tr\; ST^{-1} =\pm 2 \Leftrightarrow \lambda =
2,
\end{equation}
\begin{equation}
 \tr\; [S,T] =\pm 2
\Leftrightarrow \lambda = \pm i.
\end{equation}

J{\o}rgensen's inequality \cite{J} tells us that $G_{\lambda}$ is
not discrete if $|\lambda| < 1/2$. We refer to the circle
$|\lambda| < {\frac{1}{2}}$ as the J{\o}rgensen circle.

We can also calculate that

\begin{equation}T_{\lambda}ST_{\lambda}^{-1} =\left(
\begin{array}{cc}
1 + 2 \lambda & -4(\lambda)^2\\
1 & 1-2\lambda
\end{array}
\right)
\end{equation}
The matrices $S$ and $T_{\lambda}ST_{\lambda}^{-1}$ generate a
subgroup of $G$. We can conjugate $S$ and
$T_{\lambda}ST_{\lambda}^{-1}$ by the matrix

\begin{equation}A = \left(\begin{array}{cc}
1&0\\
{\frac{-1}{2 \lambda}} & 1
\end{array}
\right) \end{equation} so that the ordered pair of parabolics $(S,
T_{\lambda}ST_{\lambda}^{-1})$ is conjugate to the ordered pair $(S,
T_{\tilde{\lambda}})$ where ${\tilde{\lambda}}= -2 \lambda^2$. Thus
$\langle S,T_{\lambda}ST_{\lambda}^{-1} \rangle \iso
G_{\tilde{\lambda}}.$ Here $\iso$ denotes group isomorphism.
 This defines ${\tilde{\lambda}}$ for any
$\lambda$.
\begin{defn} If $\lambda \in \mathbb{C}$, we define
${\tilde{\lambda}}$ by ${\tilde{\lambda}} = -2\lambda^2$.
\end{defn}
\begin{notat}
We write $G_\l \isosup G_{\tilde{\l}}$ to indicate that $G_\l$ has a
subgroup isomorphic to $G_{\tilde{\l}}$. Similarly we write
$G_{\tilde{\l}} \isosub G_\l$. In drawing diagrams with a lattice of
subgroups, we use the standard notation
$$ G$$
$$ | $$
$$ H $$  when $H$ is a subgroup of $G$. If $G$ contains a subgroup isomorphic to
$H$,
we will write
$$ G$$
$$\; \;|\iso $$
$$ H $$
\end{notat}
We let $f(z) = -2z^2$. We will obtain information about the
structure of two-parabolic space by iterating $f$, a standard
procedure in the study of discrete group. Here, we iterate the
classical T-Schottky boundary and find potential  parabolic dust as
well as non-free points.

\section{Schottky definitions and prior results} \label{section:schottprior}
We are
interested in {\sl $\T$Schottky groups}, which are the same as
Schottky groups except that {\sl arbitrary} tangencies are allowed
among the Schottky circles and Jordan curves. That is, tangencies at
points that are parabolic fixed points and at non-parabolic fixed
points and between circles that are paired or not paired by the
marked Schottky generators are allowed.

\begin{figure}

\includegraphics[height=60mm]{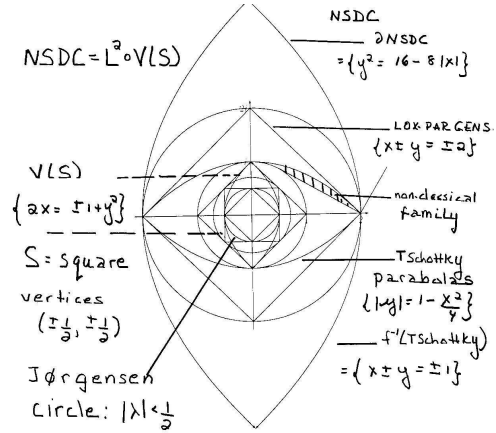}
\caption{{\bf Boundary Regions} \label{fig:rst} Each point $\lambda
\in \mathbb{C}$ corresponds to a two-generator group. The shaded
region shows $\PNC$ the one parameter family of non-classical
$\T$Schottky groups in the first quadrant. The exterior to the outer
parabolas are  the non-separating disjoint circle groups (NSDC
groups), the middle parabolas are the boundary of the classical
groups. The boundary of the Riley slice includes $\pm 2$ and $\pm i$
but otherwise lies interior to the Schottky parabolas.  Points
inside the J{\o}rgensen circle ($|\lambda| < {\frac{1}{2}}$) are
non-discrete groups. Between the Riley  slice and the J{\o}rgensen
circle are  additional non-classical groups together with degenerate
groups, isolated discrete groups and
non-discrete groups.} 
\end{figure}

The picture of the Riley slice is well known. The J{\o}rgensen
circle and the classical T-Schottky parabolas  are depicted in
figure \ref{fig:rst}.
Also depicted is a region in the plane where
all group are non-classical T-Schottky groups. Groups in the this
region which is denoted by $\PNC$  form a one complex parameter
family and are called {\sl the  non-classical family}. There are, of
course, non-classical groups outside this family. Theorem
\ref{theorem:sum} summarizes facts about these regions and their
boundaries. In this paper we iterate and apply this theorem to
two-parabolic space. First we define
\begin{defn} Let $\P = \{z=x+iy | -2 \le x \le 2 \mbox{ and } \pm y =
1-{\frac{x^2}{4}}\}$. $\P$  consists of portions of each of two
intersecting parabolas.  We refer to $\P$ as the as the {\sl
Schottky parabolas} and each portion as  a {\sl Schottky parabola}.
\end{defn}

\begin{thm} \label{theorem:sum} \cite{GW}  Let $G=
\langle A, B \rangle$ with $\Tr A =2$ and $\Tr B =2$ and $\Tr [A,B]
-2 = 4 \lambda^2$ where $\lambda \in \mathbb{C}-\{0\}$. Then

\vskip .2mm
\begin{enumerate}
\item $G$ is discrete and classical $\T$Schottky $ \Leftrightarrow \lambda = x +
iy$  and \\ $  |y| \ge 1-{\frac{x^2}{4}}.$

 \item $\l$ is on the boundary of classical $\T$Schottky space precisely when
 $\l= x+iy, -2 \le x \le 2, \mbox{ and
}   |y| = 1-{\frac{x^2}{4}}$.  That is,  the Schottky parabolas are
the boundary of classical $\T$Schottky space.

\item  \label{item:number4} If $G$ is discrete and classical $\T$ Schottky, it has a
third conjugacy class of parabolics $\Leftrightarrow$  $\lambda =
\pm 2 \mbox{ or } \pm i$. In the first case, $A^{\pm 1}B$ is
parabolic and in the second case $[A,B]$ is parabolic.

\item \label{item:hull} Let  $K$ be the convex hull of the set consisting of
the circle $|z|=1$ and the points  $z = \pm2$. If $\l$ lies between
 one of the $\T$Schottky parabolas
and the interior of $K$, then $G_\l$ is non-classical T-Schottky.
\end{enumerate}
\end{thm}

 Let $L_K$ be the
line in $\mathbb{C}$ whose equation is $3y = -{\sqrt{3}}x +
2{\sqrt{3}}$. Then $L_K$ is tangent to the circle of radius $1$ at
the point $({\frac{1}{2}}, { {\frac{\sqrt{3}}{2}} })$. Let $K_0$
denote the segment between the point of tangency and $z=2$
(including the end points) and let $C_{K_0}$ be the arc of the unit
circle between the point $({\frac{1}{2}}, { {\frac{\sqrt{3}}{2}} })$
and $i$. Set $B_{K_0} = K_0 \cup C_{K_0}$. The the region between
the interior of $K$ and the $\T$Schottky parabola in the first
quadrant is the region bounded by $B_{K_0}$ and $|y| = 1 -
{\frac{x^2}{4}}$. We denote it together with $B_{K_0}$ and its
symmetric images across the real and imaginary axes by $\PNC$ and
call it the {\sl region of the non-classical family}. We include in
$\PNC$ the line segment $K_0$ and the segment of the unit circle
between $\pi/3$ and $\pi/2$ as well as their symmetric images about
the real and imaginary axes. Figure \ref{fig:rst} shows the portion
of $\PNC$ inside the first quadrant.

 Letting $\l = x + iy = re^{i\theta}$ where $ r \ge 0$
and $0\le \theta < 2\pi$, we obtain an equivalent formulation of
\ref{item:hull}.

\begin{cor} \label{cor:region} Let $\l = x+iy = r e^{i\theta}$ lie in the first quadrant with $|y| < 1- {\frac{x^2}{4}}$.
Then $\l$ is non-classical $\T$Schottky if  either
\begin{enumerate}
\item $r \ge 1$ and
${\frac{\pi}{3}} \le \theta \le {\frac{\pi}{2}}$ or

\item $ -{ \frac{{\sqrt{3}}x}{3} } + {\frac{2{\sqrt{3}}}{3}}  \le y $ and $ 0 \le \theta \le  {\frac{\pi}{3}}$

\end{enumerate}

 If $\l$ is the image of a point in either of the above two regions
under reflection in the $x$ or $y$ axis, then $\l$ is non-classical
$\T$Schottky.

\end{cor}

 Let  $\A_1$ is the region in corollary
\ref{cor:region}(1) and $\A_2$ the region in \ref{cor:region}(2).
 We
have by symmetry the regions $- \A_i$, $\bar{\A_i}$,  and
$-\bar{A_i}$. Then $\PNC = \cup_{i=,1,2} (\A_i \cup -\A_i \cup
\bar{A_i} \cup -\bar{A_i})$.

\section{Iteration and nth T-Schottky groups}\label{section:itn}

We are interested in what happens between the J{\o}rgensen circle
and the classical $\T$Schottky parabola. We can iterate as
follows. Given $\lambda$ with ${\frac{1}{2}} < |\lambda| < 1$,
$G_{\lambda}$ has a subgroup isomorphic to $G_{\tilde{\lambda}}$
where $\tilde{\lambda} = -2 \lambda^2$. We say that $G_{\lambda}$
is {\sl second classical T-Schottky} if $G_{\tilde{\lambda}}$ is
classical T-Schottky, but $G_{\lambda}$ is not. We let $f^n(z)$
denote the nth iterate of $f$.

\newpage
 We consider the iterative sequence of subgroups,
\begin{equation}
\label{equation:subgps} \langle S, TST^{-1} \rangle
\end{equation}
$$|$$
$$\langle S, TST^{-1}\cdot S \cdot(TST^{-1})^{-1}\rangle$$
$$|$$
$$\vdots$$
$$|$$
$$\langle S,  TST^{-1}\cdot S (TST^{-1})^{-1}\cdot S \cdot (TST^{-1}\cdot
S (TST^{-1})^{-1})^{-1} \rangle$$
That is, we define inductively, the parabolic transformations
\begin{equation} \label{equation:parit}
P_0 = T,
\end{equation}
$$P_1 = TST^{-1},$$
$$P_2 = P_1 S P_1^{-1},$$
$$\vdots$$
$$P_i = P_{i-1} S P_{i-1}^{-1}.$$

\vskip .1in
 Set  $G_0= G_\l = \langle S, T_\l$ and $G_i = \langle S, P_{i-1}\rangle$. Let
$\lambda_i$ be the parameter with   $G_i \iso G_{\lambda_i}$. We
have $f^{i}(\l) = \l_i$. Note that as we move down in the lattice of
subgroups, $\l$ moves out from the origin in $\mathbb{C}$. We define

\begin{defn} If for some integer $n \ge 1, G_{f^n(\lambda)}$ is classical T-Schottky but
$G_{f^{(n-1)}(\lambda)}$ is not, we say $G_{\lambda}$ is {\sl $n$th
classical T-Schottky}.
\end{defn}

\begin{thm} \label{theorem:tncircle} Let $f(\l) = {\tilde{\l}}$ and $(f \circ f)(\l) =
{\tilde{\tilde{\l}}}$. Then
\begin{enumerate}

\item If $|\lambda| < {\frac{1}{2}},$
then $|f^n(\lambda)| < {\frac{1}{2}}$ for all integers $n$ so that
  $G_{\tilde{\lambda}}$ is non discrete and $G_{f^n(\l)}$ is
  neither discrete nor classical T-Schottky for any $n >0$.
\item  If $|\lambda| > {\frac{\sqrt{2}}{2}}$, then
$|{\tilde{\lambda}}|
> 1$ and $|{\tilde{\tilde{\lambda}}}| > 2$
  so that  $G_{{\tilde{\tilde{\lambda}}}}$ is
classical T-Schottky. \item
 Given ${\frac{1}{2}} < |\lambda| < {\frac{\sqrt{2}}{2}}$, there is
a smallest integer $n$ such that if $f(z) = -2z^2$ and $f^n(z)$
denotes the $n$th iterate of $f$, then $G_{f^n(\lambda)}$ is
classical T-Schottky but $G_{f^{(n-1)}(\lambda)}$ is not.
\item \label{item:another} For each $\l$ with ${\frac{1}{2}} < |\lambda| < {\frac{\sqrt{2}}{2}}$,
there is a unique integer $n$ such that $G_\l$ is $n$th classical
$\T$Schottky.
\end{enumerate}
\end{thm}
\begin{proof}

The first two parts are a calculation using theorem
\ref{theorem:sum} which says that if $\l \ge 2$, $G_\l$ is
$\T$Schottky.  To see the third and fourth part for each integer $n$
with
  $n\ge -1$, define $$t_n = {\frac{1}{2^{\frac{2^n-1}{2^n}}}}.$$

We consider $C_r$, the circle centered at the origin with radius $r$
where $r=t_n$. Note that $C_{-1}$ is the circle of radius $2$ and
$C_{0}$ is the circle of radius $1$. We observe that $|f(t_n)| =
t_{n-1}$ so that $f(C_n) = C_{n-1}$. Since the Schottky parabola
lies between the circle of radius $1$ and the circle of radius $2$,
the successive inverse images of the $\T$Schottky parabolas lie
between these successive circles. These circles tessellate the
region of the complex plane outside the J{\o}rgensen circle and
inside the circle of radius two as do the successive backwards
iterates of the Schottky parabolas.
\end{proof}

\section{Parabolic dust and Freeness} \label{section:dustfree}

We now investigate the pre-images of groups that are either
parabolic dust or non-free.

\begin{thm} \label{theorem:iterates}
\begin{enumerate}
\item Let $G_\l$ be such that $f^n(\l)$ is parabolic dust for some
integer $n$,
  then either $G_\l$ is parabolic dust or
non-free. \item Let $G_\l$ be such that $f^n(\l)$ is non-free for
some integer $n$, then $G_\l$ is non-free.
\end{enumerate} \end{thm}

\begin{proof} We consider the iterative sequence of subgroups $G_i
\iso G_{f^i(\l)}$  (\ref{equation:subgps}) and parabolics $P_i$
(\ref{equation:parit}) and note that as parabolic dust $G_{f^n(\l)}$
will have three distinct conjugacy classes of parabolic subgroups.
Either the generators of these groups are not conjugate in $G_\l$,
in which case $G_\l$ is parabolic dust or two or more of them are
conjugate in $G_\l$ in which case $G_\l$ has a relation and is,
therefore, non-free. As a non-free group, $G_{f^n(\l)}$ will have a
non-trivial relation and this relation will hold in $G_\l$.
\end{proof}

An immediate corollary is
\begin{cor} \label{cor:or} Let $G_\l$ be such that
$f^n(\l) = \pm i \mbox{ or } \pm 2$ for some integer $n \ge 1$
then either $G_\l$ is parabolic dust or non-free.
\end{cor}

We will prove that such groups are actually both parabolic dust and
non-free.

\comment{\vskip .2in \noindent  We also have

\begin{cor} Let $p$ and $q$ be relatively prime integers. If $p/q$ is a rational number on the real line
corresponding to a non-free point and $f^n(\l) =p/q$, then $G_\l$ is
non-free.
\end{cor}}
\subsection{Analysis of $G_\l$ for $\l = {\frac{\pm 1 \pm i}{2}}$}
 \label{section:both}

We do an analysis of $G_{{\frac{\pm 1 \pm i}{2}}}$. Note that
$G_{\tilde{\l}}\iso G_{\pm i}$ and  $G_{\tilde{\tilde{\l}}} \iso
G_{\pm 2}$. Here it is important to distinguish between subgroups
that are isomorphic to a given group and actual subgroups, that is
between $\isosub$ and $\subset$.

\begin{thm} \label{theorem:whitehead} $G_{{\frac{\pm1 \pm i}{2}}}$ is both parabolic dust and
non-free. It contains two conjugacy classes of $\mathbb{Z} \times
\mathbb{Z}$ parabolic subgroups and thus is not free. It contains
four distinct conjugacy classes of parabolic elements and thus is
parabolic dust.

\end{thm}

\begin{proof}
By symmetry, it suffices to  treat the case $\l = {\frac{1-i}{2}}$.
\begin{enumerate}
\item  $\l= {\frac{1-i}{2}}$,  $G_{\lambda} = \langle S,
T_{\lambda} \rangle$,  $T_\l$ fixes $\infty$ and $S$ fixes $0$.

 \item $f(\l) = -2\l^2 =i = {\tilde{\l}}$,  $G_{\tilde{\l}} \iso
\langle S, T_\l S T_\l^{-1} \rangle$ \\
 Let $\tilde{T}= T_\l S
T_\l^{-1}$. $\tilde{T}$ fixes $1-i$ and $S$ fixes $0$. \\
We note that $S\tilde{T}S^{-1}{\tilde{T}}^{-1} = \left(%
\begin{array}{cc}
  3 & -4 \\
  4 & -5 \\
\end{array}%
\right)$  is parabolic with fixed point at $1$.

\item  ${\tilde{\tilde{\l}}} =2$, $G_{\tilde{\tilde{\l}}}=
G_{2}$ and $\tilde{\tilde{T}} = {\tilde{T}} S {\tilde{T}}^{-1} =
\left(%
\begin{array}{cc}
  -1 & 4 \\
  -1 & 3 \\
\end{array}
\right)$ has its fixed point at $2$.

\item We look at $S^{-1}{\tilde{\tilde{T}}} = \left(%
\begin{array}{cc}
  -1 & 4 \\
  0 & -1 \\
\end{array}%
\right)$. Its fixed point is at $\infty$.

\item Thus $T_{\lambda}$ and $S^{-1}{\tilde{\tilde{T}}}$ are parabolics that share the
same fixed point, but are not conjugate in the M\"obius group.
Therefore, they must commute \begin{equation}
\label{equation:commute} T_\l S^{-1}{\tilde{\tilde{T}}} =
S^{-1}{\tilde{\tilde{T}}}T_\l.
\end{equation} and, therefore, $G_{{\frac{1-i}{2}}}$ is not free
and contains the  $\mathbb{Z} \times \mathbb{Z}$ subgroup generated
by the commuting elements.
\end{enumerate}

\noindent We note that a similar analysis interchanging $S$ and
$T_\l$ shows that $S$ and $T_\l^{-1}(S T_\l^{-1} S^{-1} T_\l \cdot
(ST_\l S^{-1})^{-1}$ share a fixed point, namely $0$. Since $G_\l$
is non-elementary, $S$ and $T_\l$ are not conjugate, none of the
four parabolics can be conjugate.
\end{proof}

Looking at further pre-images of $\pm i$, we have
\begin{cor} \label{cor:generalwhite} If $f^n(\l) = \pm i$ for some $n>1$, then $G_\l$ is
simultaneously non-free and parabolic dust and contains at least two
distinct $\mathbb{Z} \times \mathbb{Z}$ subgroups.
\end{cor}
\begin{proof} If $f^n(\l) = \pm i$, then $f^{(n-1)}(\l) = {\frac{\pm
1 \pm i}{2}}$ so that  $G_\l \isosup G_{{\frac{\pm 1 \pm i}{2}}}$
and inherits its two non-conjugate $\mathbb{Z} \times \mathbb{Z}$
subgroups.
\end{proof}
One could further analyze the pattern  of $\mathbb{Z} \times
\mathbb{Z}$ subgroups in the pre-images of $\pm i$. We expect the
number of non-conjugate such groups to grow.

\begin{figure}
\includegraphics[height=60mm]{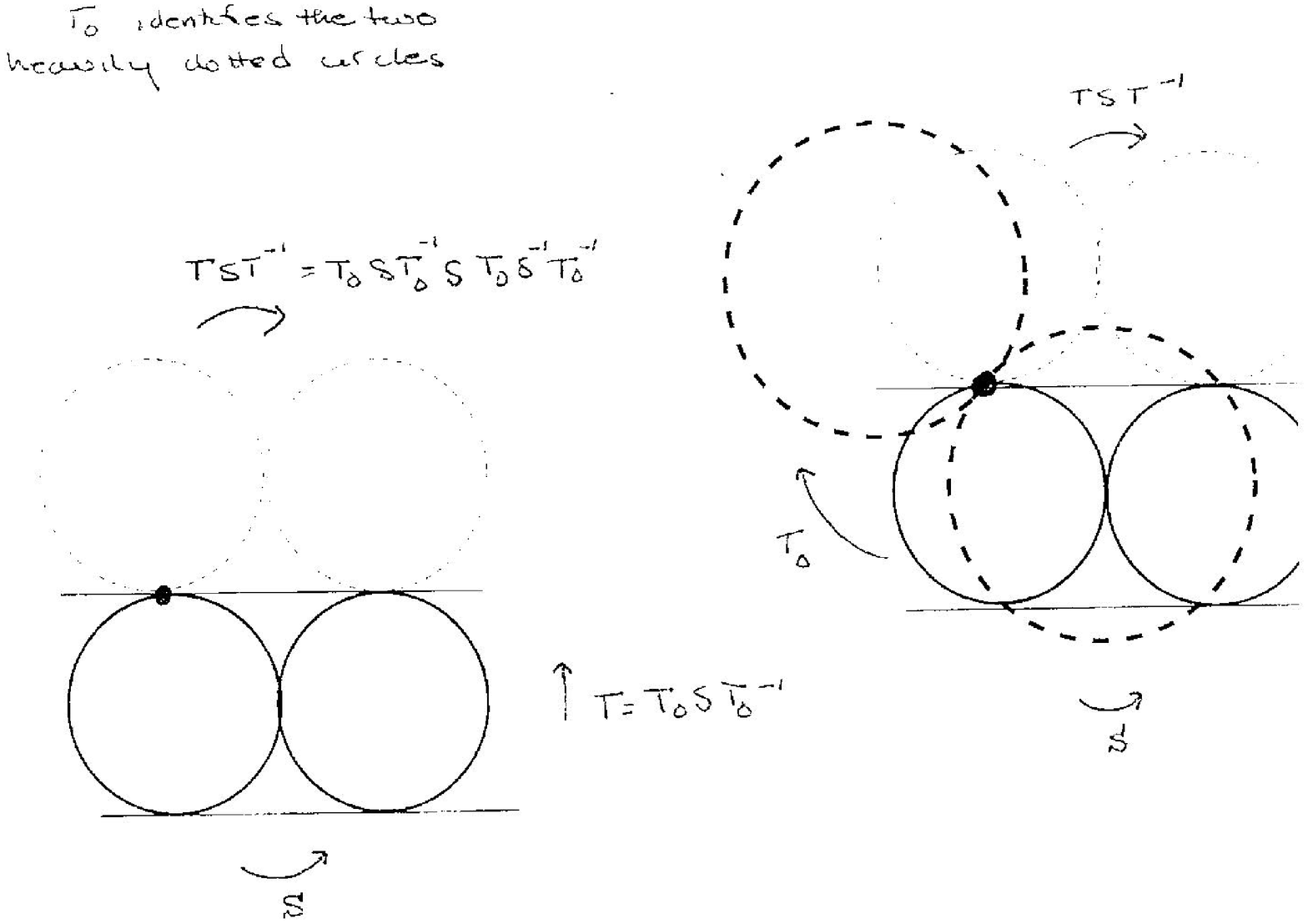}
  \caption{We have conjugated the group
$G_{{\frac{1-i}{2}}}$
\comment{
\begin{array}{cc}
  -1 & 0 \\
  -{\frac{-1+i}{2}} & -1 \\
\end{array}%
\right)$} so that $S$ identifies $C = \{ |z+ 1| =1\}$  and
$C'=\{|z-1|=1$, $T$ identifies $L =\{ y= -i\}$  and $L'= \{ y=i \},
$
 $S(z) = {\frac{z}{z+1}}$, $T(z) = z + 2i$,
 $TST^{-1}(z) = {\frac{(1+2i)z + 4}{z - 2i +1}}$,
  $T_0(z) = {\frac{2z + 1-i}{ {\frac{-1+i}{2}} z}}$, $T= T_0ST_0^{-1}$,
 $G_{\frac{-1+i}{2}}
  \iso \langle S, T_0 \rangle \supset   G_i = \langle S, T \rangle
 \supset G_2 = \langle S, TST^{-1} \rangle$ Then $T_0$ has the same fixed point
 as $[S^{-1},T]$,
    \comment{
\left(%
\begin{array}{cc}
  1 + 2i & 4 \\
  -2i & -3+ 2i \\
\end{array}%
\right)$}
 namely $-(1+i)$.} \label{fig:subs}
\end{figure}

Figure \ref{fig:subs}
illustrates the action of (a conjugate of)
$G_{{\frac{1-i}{2}}}$ and its subgroups in ${\hat{\mathbb{C}}}$,
showing the parabolic fixed points and the isometric circles of the
transformations. We have conjugated so that $\infty$ is no longer a
parabolic fixed point.

\begin{rem} The $G_{\frac{\pm 1 \pm i}{2}}$ correspond to the Whitehead link and
have finite volume quotients. Since at each step the groups are of
infinite index, the corresponding covering would have infinite
degree, so the $G_\l$ with $f^k(\l) = {\frac{\pm 1 \pm i}{2}}$ will
have zero volume. That is, the limit sets will be
${\hat{\mathbb{C}}}.$
\end{rem}

\subsection{Additional Non-free points}

 We note that Beardon, Bamberg, Ignatov, and Lyndon-Ullman
\cite{Beard, BPell, BAM, Ignatov, IgnatovII,  LU} among others have
studied non-free points on the real axis. We combine their results
along with theorem \ref{theorem:iterates} to obtain complex non-free
points. \newpage
\begin{thm} \label{theorem:nonfree} Let $\delta\ne 0 \in \mathbb{C}$. Then $G_{\delta}$ is not free if for some integer $k \ge 0$
$f^k(\delta) = \l$ and  one of the following holds:
\begin{enumerate}

\item 
$\l \in \{ {\frac{p^2}{2(p^2+1)^2}} \; |
\; p=1,2, ...\} \cup \{{\frac{1}{2p^2}}\;|\; p=2,...\}$

\item 
  $\l = 4\sin^2({\frac{p\pi}{q}})$ where $p$ and $q$ are relatively
prime integers.

\item 
$\l =
{\frac{p^2}{2q^2}}$ where  $p$ and $q$ satisfy Pell's equation,
$1=q^2-Np^2$,  for some positive square free integer  $N$.
\item $\l
= {\frac{p_n^2}{2q_n^2}}$, $\l$ where $N$ is a positive square free
integer and the convergents of ${\frac{1}{\sqrt{N}}}$ are $p_n/q_n$.

\item 
$\l \in \mathcal{B}$  where $${\mathcal{B}} =
 \{
 {\frac{1}{2}},
  1,
 {\frac{3}{2}},
{\frac{1}{2}}e^{{\frac{i \pi}{3}}},
 {\frac{1+ {\sqrt{13}}}{4}},
 {\frac{5+ \sqrt{5}}{4}}, {\frac{i}{\sqrt{2}}},{\frac{1}{\sqrt{2}}}, {\frac{ 9 }{ 50   }},
  {\frac{8 }{ 25}},
{\frac{25}{72}}, {\frac{8}{81}}, {\frac{25}{162}}, {\frac{25 }{98
}}, {\frac{9 }{8 }}, {\frac{ 8  }{ 9   }},
 {\frac{25   }{32 }},
{\frac{25 }{18}}  \}.$$
\item 
 $\l \in
  {\frac{\l_0}{n^2}}
  $
where
$G_{\l_0}$ is not free.

 \item 
  $\l \in  \{
  {\frac{1}{2}}({\frac{1}{n}})^2,
  {\frac{1}{2}}({\frac{2}{n}})^2,
  {\frac{1}{2}}({\frac{3}{n}})^2,...,
{\frac{1}{2}}({\frac{8}{n}})^2,  \; \; n\ne 0 \in \mathbb{Z} \}$
\item  $\l \in  \{{\frac{(m+n)^2 }{2m^2n^2 }}\; |\;  m\ne0 ,n\ne0 \in
\mathbb{Z} \}$
\end{enumerate}
\end{thm}

\begin{proof} Apply theorem  \ref{theorem:iterates} to the results of Bamberg, Beardon, Ignatov and
Lyndon-Ullman. Translate from the $\mu$ parameter to the $\l$
parameter via $2\l=\mu^2$ where necessary.
\end{proof}

\subsection{Pre-images of cusps on the boundary of the Riley slice}
In \cite{KS} and \cite{Wright, Wr} a family of words $W_{p/q}$
indexed by the rationals,  $p/q$,  are studied. These are termed
Farey words. For a discrete group, words in the group correspond to
curves on the corresponding surface. In particular pinching the
curves corresponding to Farey words until $\tr W_{p/q} = -2$ gives
cusps on the boundary of the Riley slice, that is,  parameters where
there are three parabolics. We let $\l_{p/q}$ be the value of $\l$
for which $\tr W_{p/q} = -2$. Apply theorem \ref{theorem:iterates}
to see that
\begin{cor}
Let $G_\l$ be such that $f^n(\l) = \l_{p/q}$ for some integer $n$
and some rational number $p/q$. Then either $G_\l$ is parabolic dust
or non-free.
\end{cor}

\subsection{Questions}
Ian Agol has asked the question about the values of $|\l_p/q|$.
David Wright has computed $|\-p/q|$ for some values of $p/q$. It
would be nice to understand where these points lie in the in the
tessellation of the space by the circles $C_{t_n}$ that arose in the
proof of theorem \ref{theorem:tncircle}.

Can we decide whether the group has a third parabolic conjugacy
class or whether it is non-free? Bamberg \cite{BAM} has a procedure
for determining whether $\lambda$ represents a non-free group. It is
not an algorithm, that is, it does not necessarily stop, but if one
has a bound on the length of the word, which we do, and a bound on
the exponents of the word in the original generators, his procedure
may be useful.

\section{Tessellation of $\mathbb{C}$ by  the  backwards iterates of the Schottky
parabolas} \label{section:tess}

With some more terminology we can say a number of useful things
about the backwards iterates of the Schottky parabolas and other
boundaries and regions.

\begin{defn}
Let $D_{\diamond}$ be the Euclidean square in $\hat{\mathbb{C}}$
with vertices $\pm i$ and $\pm 1$.
\end{defn}

\subsection{The pre-image of the Schottky parabolas }

To investigate the  pre-images of  $\P$, we first look at what $f$
does to lines.  We will see that the function $f$ maps horizontal
and vertical lines to parabolas with vertices on the $x$-axis. It
also maps all other straight lines to parabolas.

\subsubsection{Images of lines}
 We want to investigate the image of the straight line
$x+y=c$ under the map $f(z) = -2z^2$ where $c$ is a real number and
$z=x+iy$. We set $f(z) = \tilde{z} = \tilde{x} + i\tilde{y}$.
\begin{lem} \label{lem:prelines} The image of the line $x+y =c$ under $f$, is the parabola $\tilde{y}= -c^2 + {\frac{(\tilde{x})^2}{4c^2}}.$

\end{lem}

\begin{proof}  If $f(z) = \tilde{z}= \tilde{x} + i\tilde{y}$, $\tilde{x} =
-2(x^2-y^2)$ and $\tilde{y}= -4xy$. Calculate that if $x+y=c$,
$y=c-x$,
$$y^2=c^2-2cx+x^2$$
$$x^2-y^2 = 2cx-c^2$$
$$\tilde{x} = -2(x^2-y^2)=2c^2-4cx$$
$$(\tilde{x})^2 = 4c^4-16c^x + 16c^2x^2$$
and
$${\frac{(\tilde{x})^2}{4}} =c^4-4c^3x+4c^2x^2.$$

Also calculate that
$$\tilde{y} = -4xy = -4x(c-x)=-4xc + 4x^2$$

$$ c^2 \cdot \tilde{y} = -4xc^3 + 4x^2c^2$$

$$ c^2 \cdot \tilde{y} + c^4 = c^4 -4xc^3 + 4x^2c^2=
{\frac{(\tilde{x})^2}{4}}$$

$$ c^2 \cdot \tilde{y} = -c^4 + {\frac{(\tilde{x})^2}{4}}$$

$$\tilde{y}= -c^2 + {\frac{(\tilde{x})^2}{4c^2}}$$
\end{proof}

Whence,
\begin{cor} \label{cor:Sparab} If $c= \pm 1$, $\tilde{y}= -1 +
{\frac{(\tilde{x})^2}{4}}.$
\end{cor}

\begin{cor} The diamond, $D_{\diamond}$,  is the preimage of the Schottky
parabolas.
\end{cor}

\begin{proof} Use corollary \ref{cor:Sparab}
A similar calculation to the one done in proving lemma
\ref{lem:prelines} shows that if $x-y = \pm 1$, $\tilde{y}= 1 -
{\frac{(\tilde{x})^2}{4}}.$
\end{proof}

Similar calculations yield

\begin{lem} The line $x_0=c$ is mapped to the parabola  $x = -2c^2
+{\frac{y^2}{8c^2}}$.

The line $y_0=c$ is mapped to the parabola $x = 2c^2
-{\frac{y^2}{8c^2}}$.

\end{lem}

\comment{
\begin{proof}
If $z_0 = c+ iy_0$ and $f(z_0) = -2z_0^2 = z = x+ iy$, then $x =
-2c^2 + 2y_0^2$ and $y= -4cy_0$. Thus $x =-2c^2 -
{\frac{y^2}{8c^2}}$. The calculation for $y_0=c$ is similar.
\end{proof}}
In particular,
\begin{cor}
The function $f$ maps the square with vertices at $(\pm c,\pm c)$ to
the parabolas $\tilde{x} = \pm 2c^2 + {\frac{\mp
\tilde{y}^2}{8c^2}}$
\end{cor}
and \begin{cor} \label{corollary:S} The function $f$ maps the square
with vertices at $(\pm {\frac{1}{2}}, \pm {\frac{1}{2}})$ to the
region bounded by the parabolas ${\tilde{x}} = {\frac{\pm 1}{2}} +
{\frac{\mp 1}{2}}(\tilde{y})^2.$
\end{cor}

 \comment{This shows that
if $D_{\diamond}$ is the region bounded by the line $y=x+1$,
$y=-x+1$, $y=-x-1$ and $y=x-1$. That is, $D_{\diamond}$ is the
square with vertices $\pm 1$ and $\pm i$, (sides of length
$\sqrt{2}$),  then $f(D_{\diamond})$ is "the T-Schottky parabolas".
by abuse of language we call the boundary}

\subsubsection{Pre-images of lines}

In order to identify the images and pre-images of further regions
such as $\PNC$ under iteration, we consider pre-images of the line
$y = mx +b$ where $m \ne 0$. In looking at the regions, we will
usually only be interested in segments of lines or parabolas or
hyperbolas as we were when we identified the Schottky parabolas as
segments of two intersecting parabolas. We adopt the convention the
convention that when we say the preimage or image is a line or a
parabola or a branch of a hyperbola, we mean the appropriate segment
which should be clear from the context.

\begin{lem} Let $L = \{z = x+iy\;|\; y=mx+b\}$ be a straight line with $m \ne 0$ and $b \ne 0$,
then $f^{-1}(L)$ is the hyperbola whose equation is

$$ {\frac{(x-{\frac{y}{2m}})^2}{({\frac{1}{2m}}+2m)b}} - {\frac{y^2}{2mb}} = 1.$$

\end{lem}
\begin{proof} Let $f(x_0 + iy_0) = z = x+ iy$, then $x =
-2(x_0^2-y_0^2)$ and $y= -4x_0y_0$. Substitute into $y=mx+b$,
simplify, and then replacing
$x_0$ by $x$ and $y_0$ by $y$.

\end{proof}

 Recall that the line $L_K$  is the line with the equation $y=
-{\frac{\sqrt{3}}{3}}x + 2{\frac{\sqrt{3}}{3}}$. We call this the
$K$-line.  The point of tangency of the $K$-line and the unit circle
is $({\frac{1}{2}}, {\frac{\sqrt{3}}{2}})$.

\begin{ex} \label{lem:Kline} {\rm
The preimage of the segment of the $K$-line between $x =
{\frac{1}{2}}$ and $x=2$ is the portion of the hyperbola $4(x_0 +
{\frac{\sqrt{3}}{2}})^2 - 7y_0^2 = 3$ between $i
{\frac{\sqrt{2}}{2}} e^{i{\frac{\pi}{6}}} $ and $-i $.}
\end{ex}

\begin{ex}\rm{Let $L_{\theta}$ denote the line through the origin making a
counterclockwise angle of  $\theta$ with the positive real axis.
Then $f(L_{\theta}) = L_{2\theta}$.}
\end{ex}

\begin{ex} \rm{Recall that  $\A_1$ is the region in the first quadrant bounded by
the unit circle, the Schottky parabola, and $L_{{\frac{\pi}{3}}}$.
Thus $f^{-1}(\A_1)$ will be bounded by
 the circle centered  at the
origin with radius ${\frac{\sqrt{2}}{2}}$, the line $x+y=1$ and
$L_{{\frac{\pi}{6}}}$.}
\end{ex}

\subsection{Second pre-image of the Schottky parabolas}
We have seen that  $f(\{x+iy \;| \;\;x \pm y = \pm 1\})$ is the
Schottky parabolas. We want to describe $x_0 + iy_0$, where
$\{f(z_0)\} = \{ x+iy\;| \;\;x \pm y = \pm1\}$.

We note that $x = -2(x_0^2 -y_0^2)$ and $y = -4x_0y_0$. Thus
$$-2(x_0^2-y_0^2) \pm (-4x_0y_0) = \pm1.$$ In particular,
$${\frac{(x_0 \mp y_0)^2}{\pm {\frac{1}{2}}}} \mp
{ \frac{(y_0)^2}{ \mp {\frac{1}{4}}} }=1.$$
This proves 
\begin{cor} \label{cor:hyps}
The second inverse image of the Schottky parabolas consists of
portions of eight hyperbolas. The hyperbolas intersect in pairs on
points on the circle $C_{\frac{\sqrt{2}}{2}}= f^{-2}(C_0)$ and the
portions of the hyperbolas between sucessive points of intesection
gives the boundary of second classical $\T$Schottky space.
\end{cor}

\begin{figure}
\includegraphics[height=60mm]{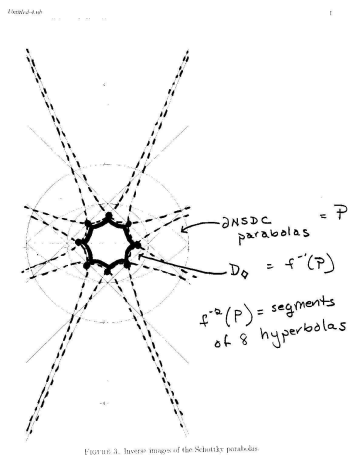}
\caption{Inverse images of the Schottky parabolas.}
\label{fig:mathematica} \end{figure}

Figure \ref{fig:mathematica}
 illustrates this.  More
generally, if $n \ge 2$, one expects $f^{-n}(\P)$, the boundary of
the $n$th classical space, to be the union of portions of branches
of $2^{n+1}$ hyperbolas. These $2^{n+1}$ branches of hyperbolas
intersect in pairs with the points of intersection lying on the
circle $f^{-n}(C_0)$ which has radius $t_{n-1} =
{\frac{1}{2^{\frac{2^{n-1}-1}{2^{n-1}}} }}$ and the boundary
consists of the portions of the hyperbolas between successive
intersection points. The $\{f^{-n}(\P), n \in \mathbb{Z} \}$
tessellate the complex plane.

\begin{ex} {\it Loxodromic-Parabolic T-Schottky generators.}
{\rm Any given pair of generators of a classical Schottky or
$\T$Schottky group may or may not be generators for an actual
Schottky configuration, that is, the set of side pairings for two
pairs of Schottky circles. If they are we say that they are {\sl
Schottky generators for a marked Schottky configuration}. In
\cite{GW} it was shown that if a $G_\l$ was a $\T$Schottky group, it
 had a pair of parabolic Schottky generators for a marked Schottky
configuration. It was also shown that $G_\l$ had an additional
loxodromic-parabolic pair of Schottky generators  for a marked
Schottky configuration precisely when $\l = x + iy$ satisfied
 $x \pm y \ge \pm 2$ \cite{GW}. In particular if $D_{2\diamond}$ is the
square with this boundary its pre-image under $f$ consists of
segments of the eight hyperbolas $\pm (x_0 \mp y_0)^2 \mp
{\frac{y_0^2}{\mp {\frac{1}{2}} }} = \pm 1$ inside the unit circle
that intersect at the  $\pm 1$, $\pm i$, and ${\frac{\pm 1\pm
i}{2}}$. See Figure \ref{fig:rst}.}
\end{ex}

\section{Lateral Iteration and Vertical iteration} \label{section:lateral}

We think of the iteration of subgroups given by
\ref{equation:subgps} as being {\sl vertical iteration}. The
successive subgroups are moving down in the diagram of the lattice.
We can also perform what we call a {\sl lateral iteration}. That is,
for each integer $n \ne 0$, $T_{n\l} = T_\l^n$ and  we consider the
subgroup $\langle S, T_{n\l} \rangle$ of $G_{\l}$. Now $\langle S,
T_\l^n \rangle = \langle S, T_{n\l}\rangle = G_{n\l}$.

We can also perform a single  vertical iteration or a series of
vertical iterations  on  $G_{n\l}$. A single vertical iteration
replaces $G_{n\l}$ by $G_{-2{n\l}^2} = G_{-8\l^2}$ and  $m >0$,
vertical iterations replaces $G_{n\l}$ by $G_{f^m(n\l)}$.

Since vertical iteration replaces $\l$ by $-2\l^2$ and lateral
iteration replaces $\l$ by $n\l$, we think of $V$ moving groups down
in the lattice of subgroups more quickly than $L$. This motivates
the names.

 If we let
$V(G_\l) = G_{\tilde{\l}}$ and $L^n(G_\l) = G_{n\l}$, then we can
describe successive vertical and horizontal iterations by a sequence
of $V^{m_1}L^{n_1}V^{m_2} \cdots L^{n_t}$ where $m_1,...m_t$ and
$n_1,...,n_t$ are positive integers and $V^m$ denotes $m$ successive
vertical iterations so that $V^m(G_\l) = G_{f^m(\l)}$.

 Combining horizontal and vertical iteration allows one to move
more fully around in the parameter space. We observe that
\begin{prop} \label{prop:lv}
\begin{enumerate}
\item No applications of vertical iteration moves $\l$ out of the
J{\o}rgensen circle. \item There always is a lateral move that takes
a given $\l$ out of the J{\o}rgensen circle. That is,  for each $\l
\ne 0$ in the J{\o}rgensen circle, there is a smallest $m$ such that
$L^m(G_\l)= G_{m\l}$ lies outside the J{\o}rgensen circle. That is,
there is an integer $m$ such that $|m\l| > {\frac{1}{2}}$.
\item $L$  and
$V$ do not commute in general. That is if $n \ne 1$,  \\$L^n\circ V
(G_\l) \iso G_{-2n\l^2} \ne G_{-2n^2\l^2} \iso V \circ L^n (G_\l)$.
\item  However, there are some relations among the $L$ and $V$. For example
 \begin{equation} V\circ L^m(G_\l) = L^n \circ V(G_\l)
 \Leftrightarrow n = m^2.
  \end{equation}
\end{enumerate}
\end{prop}

One can use lateral iteration and obtain more formulas for parabolic
dust and non-free groups.

\begin{ex}{\it The NSDC parabolas under combined Vertical and Lateral iteration:}
{\rm In \cite{GW} non-separating circle groups, NSDC groups for
short, were studied. It was found that the boundary of NSDC space is
given by the portions of the two parabolas $y^2 = 16\pm 8x $ for $-2
\le x \le 2$. Let $S$ be the square with vertices at $(\pm
{\frac{1}{2}}, \pm {\frac{1}{2}})$. Then $V(S)$ consists of the
parabolas $2x = \pm 1 + y^2$. Apply $L^{-2}$ (e.g. $ z \mapsto -4z$)
to these  parabolas to obtain the NSDC parabolas setting as setting
$x= -{\frac{x_0}{4}}$ and $y = -{\frac{y_0}{4}}$  yields $y_0^2 = 16
+ 8|x_0|$. That is, $L^{-2}\circ V(S) = \partial NSDC.$ See Figure
\ref{fig:rst}.}
\end{ex}

\section{The action of $S$ and $T$ on ${\hat{\mathbb{C}}}$.}
\label{section:configs}

We have described the various types of groups by the regions in
which $\l$ lies. We can equally well describe these groups by the
configuration of Schottky curves imposed upon a group with a given
type of $\l$.

\begin{figure}
\includegraphics[height=60mm]{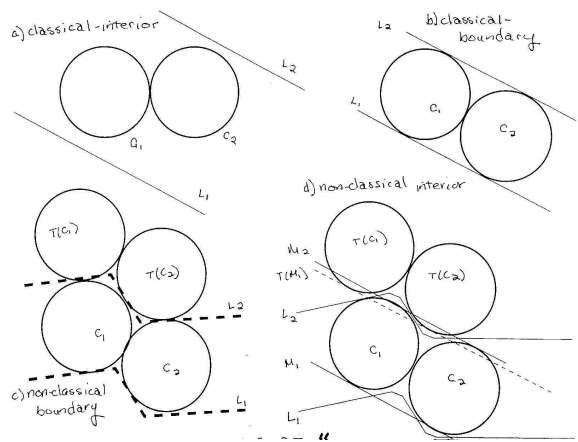}
 \caption{This figure illustrates the domains: a) the interior classical $\T$Schottky
 configuration, b) the boundary classical $\T$Schottky configuration,
  c) the interior non-classical configuration in the $\PNC$ family
  and d) the non-classical boundary configuration for the $\PNC$
  family. }
\label{fig:config}
\end{figure}

If a rank two Schottky group has generators $g_1$ and $g_2$ and
$g_1(C_1)=C_2$ and $g_2(T_1) = T_2$ where $C_1,C_2,T_1, \mbox{ and }
T_2$ are Jordan curves whose interiors are pairwise disjoint, then
the intersection of the exterior of these curves is called the {\sl
Schottky domain} and the set of circles together with the side
pairings is called a {\sl Schottky configuration}. If the Jordan
curves are circles, we have a {\sl classical Schottky
configuration}. Not every set of generators for a Schottky group are
the side pairings for a Schottky configuration. If $g_1$ and $g_2$
are the side pairings for a Schottky configuration, we call them the
{\sl Schottky generators}. In the case of $G_\l$ we normalized so
that $C_1$ and $C_2$ paired by $S$ were are tangent at $0$ with
equal radii, $R$, and $T_1$ and $T_2$ paired by $T$ were tangent at
$\infty$. We summarize the geometric configurations obtained for
various regions. See Figure \ref{fig:config}.

\begin{thm} A $G_\l$ for $\l$ in the given region
always has a Schottky configuration that can always be constructed
as indicated.
\vskip .2in

 \begin{description}
    \item[$\l$ is an interior $\T$Schottky point] One can always find a classical configuration where the
   where the  circles $C_1$ and $C_2$ lie interior to $T_1$ and $T_2$ and there are no additional
 tangencies between the circles. The configuration has a total of
 two points of tangency.
    \item[$\l$ is $\partial \T$Schottky] Every classical Schottky configuration will have
    six tangencies between
    the four circles except when $\l=\pm 2$ when there are only fours tangencies.

    \item[$\l$ lies on the lower boundary of $\PNC$, on the $K$-line] (or on the arc of
    the unit circle).  The curves $C_1,C_2,T(C_1),T(C_2)$ will give a classical Schottky
    configuration for $\langle S, TST^{-1}
    \rangle$.  Here $C_1$ and $C_2$ are round circles. $T(C_1)$ is tangent to $C_1$
    at a point $p_1$,   $T(C_2)$ is tangent to $C_2$ at a point $p_2$,  and $C_2$ is tangent
    to $T(C_1)$ at a point $p$. Let $M_1$ be the tangent line to $C_1$ through
    $p_1$,  $M_2$ the tangent line to $C_2$ through $p_2$;  and
    $M_3$ the tangent line to $T(C_1)$ and $C_2$ through $p$. Let $L_2$ be the Jordan curve made up of three arcs:
    the segment of $M_1$ between $\infty$ and $M_1\cap M_3$,
    the segment of $M_3$ between $M_1 \cap M_3$ and $M_3 \cap M_2$, and
     the segment $M_2$ between $M_3 \cap M_2$ and $\infty$. Set $L_1 = T^{-1}(L_2)$.
    Then the transformation $T$ pairs the Jordan curves $L_1$ and
    $L_2$.
\item[$\l$ interior to $\PNC$]

Let $M_1$ and $M_2$ be parallel lines with each tangent to both
$C_1$ and $C_2$.  We assume that $T(M_1)$ lies between  $M_1$ and
$M_2$. \comment{ the line connecting the centers of $C_1$ and $C_2$
and $M_2$. Let $P$ be the perpendicular bisector of the line
    connecting the centers of $C_1$ and $C_2$. Let $T(M_1)$
    intersect $P$ at a point $q_1$ and let $M_2$ intersect $P$ at $q_2$.
    Let $q$ be any point on $P$ between $q_1$ and $q_2$ and let $q_0$ be the point on $P$
    that is $R$ units from $q$ where $R$ is the radius of $C_1$.
    Let $D_1$ be the circle of radius $R$ centered at $q_0$ and
    $D_2= S(D_1)$. Then taking segments of the perpendicular
    bisectors of the lines between the centers of $C_1$ and $D_1$,
    etc. will work.}Let $P$ be the perpendicular bisector of the line connecting the
centers of $C_1$ and $C_2$. Let $M$ be a line parallel to $M_1$ that
intersects $P$ at $q_0$, a point between the point of tangency of
the $C_1$ and $C_2$ and $M_2 \cap P$. Let $q$ be the point on $P$ a
distance $R$ from $q_0$. Let $2\l$ be the vector from $0$ to $q$ and
let $T(z) = z + 2\l$. Let $L_2$ be the Jordan curve comprised of
arcs along the perpendicular bisector of the centers of $C_1$ and
$T(C_1)$, along $M$, along the perpendicular bisector of the centers
of $T(C_1)$ and $C_2$, along the line parallel to $M$ equidistant
from $M_2$, and along the perpendicular bisector of the line segment
connecting the centers of $C_2$ and $T(C_2)$. Let $L_1=T^{-1}(L_2)$.
Then $L_1$, $L_2= T(L_1)$, $C_1$ and $C_2= S(C_1)$ give a Schottky
domain for $G_\l$ and $\l$ will be interior to $\PNC$.
 \end{description}
\end{thm}
\begin{proof}
The description for the first three regions come from \cite{GW}. The
description for the last follows by ruling out the configurations
for first three cases.
\end{proof}

\section{Towards a more general theory} \label{section:conclude}

If one looks in the bigger space of representations modulo conjugacy
for two generator groups, what we have termed {\it lateral} is
really {\sl semi-lateral iteration}, $SL$. It is lateral only in
comparison with $V$.  A general Nielsen transformation  sends
generators $(A,B)$ to generators $(A, A^{\pm 1}B)$,  does not change
the group, and does not change the properties of discreteness or
non-discreteness or free-ness. This is true lateral iteration. It
does change the triple of parameters one associates to a two
generator group. We note that in the case of a discrete group, such
a Nielsen transformation, $N$, corresponds to an element of the
mapping-class group and thus has a geometric interpretation as a
Dehn twist.

In our case, the Nieslen transformation moves us out of the
two-parabolic parameter space. But when we move to a more general
theory,  it might be that certain sequences of moves in the larger
parameter space bring us back to the two-parabolic space. We would
like to develop a theory of motions around the representation space
using the three types, $V$, $SL$ and $N$. This would be related to
the various word families studied in connection with representations
and volumes of hyperbolic manifolds and orbifolds.

\section{Acknowledgements} The author thanks David Wright for a number of helpful
e-mail exchanges. The author also thanks Yair Minsky and the Yale
Mathematics Department for their hospitality and support while some
of this work was carried out.

\vskip .07in

Mathematics Department, Rutgers University, Newark, NJ 07102

e-mail: gilman@andromeda.rutgers.edu







\end{document}